\documentclass{amsart}
\newtheorem{theorem}{Theorem}[section]

\theoremstyle{definition}

\newtheorem{example}[theorem]{Example}

\theoremstyle{remark}
\newtheorem{remark}[theorem]{Remark}

\numberwithin{equation}{section}

\newcommand{\R}{{\mathbb R}}
\newcommand{\C}{{\mathbb C}}
\newcommand{\N}{{\mathbb N}}

\newcommand{\cP}{{\mathcal P}}

\newcommand{\sgn}{\operatorname{sgn}}

\renewcommand{\epsilon}{\varepsilon}
\newcommand{\cal}[1]{{\mathcal{#1}}}

\newcommand{\supp}{\operatorname{supp}}

\begin{document}

\title[Dynamics of Prions II]
{Analysis of a Model for the Dynamics of Prions II}

\author{Hans Engler}
\address{Department  of Mathematics, Box 571233, Georgetown University, Washington D.C. 20057-1233, USA}
\email{engler@georgetown.edu}
\author{Jan Pr\"uss}
\address{Martin-Luther-Universit\"at Halle-Wittenberg,
Fachbereich Mathematik und Informatik, Institut f\"ur Analysis,
Theodor-Lieser-Str. 5, D-06120 Halle, Germany}
\email{jan.pruess@mathematik.uni-halle.de}
\author{Glenn F.~Webb}
\address{Department of Mathematics,
Vanderbilt University, Nashville TN 37212, USA}
\email{glenn.f.webb@math.vanderbilt.edu}

\vspace{1cm}

\maketitle

\begin{center}
{\bf Abstract}
\end{center}

\noindent

{\footnotesize  A new mathematical model for the dynamics of prion
proliferation involving an ordinary differential equation coupled
with a  partial integro-differential equation is analyzed,
continuing the work in \cite{PPWZ05}. We show the well-posedness
of this problem in its natural phase space $Z_+:=\R_+\times
L_1^+((x_0,\infty);xdx)$, i.e. there is a unique global semiflow
on $Z_+$ associated to the problem.

A theorem of threshold type is derived for this model which is
typical for mathematical epidemics. If a certain combination of
kinetic parameters is below or at the threshold, there is a unique
steady state, the disease-free equilibrium, which is globally
asymptotically stable in $Z_+$; above the threshold it is
unstable, and there is another unique steady state, the disease
equilibrium, which inherits that property. }


\date{March 28, 2005}

\begin{center}
{\bf Acknowledgement}
\end{center}

\noindent

{\footnotesize This paper was initiated while the second author
was visiting the Department of Mathematics, Vanderbilt University,
Nashville, Tennessee in  2003/04. He wants to express his thanks
to the department for kind hospitality and for financial support.}

\section{Introduction and Main Results}

In this paper we continue our analysis, begun in \cite{PPWZ05}, of
a recent model describing the proliferation of prions. This model
has been introduced in Greer, Pujo-Menjouet and Webb \cite{GPW04},
based on the works of Masel, Jansen and Nowak \cite{MJN99}, Nowak,
Krakauer, Klug and May \cite{NKKM98} and others. For comprehensive
explanations and discussions of the model and the relevant
biochemical  literature we refer to \cite{GPW04}. Here we only
give a very short description of the model.

Prions are proteins that are believed to be responsible for
certain diseases  like BSE and the Creutzfeld-Jacob disease. There
are two basic forms of prions of interest here, the {\em Prion
Protein Cellular} $PrP^C$ and the {\em Prion Protein Scrapie}
$PrP^{Sc}$. The single molecule proteins $PrP^C$, also called {\em
monomers} in the sequel, are protease resistent proteins which
have a cell protective function and are produced by the body,
regularly. On the other hand, the infectious prion $PrP^{Sc}$ is a
string-like {\em polymer} formed of monomeric $PrP^C$. Above a
critical chain length $x_0>0$ the polymers are more stable than
the $PrP^C$, and they can grow to chains containing thousands of
monomers. $PrP^{Sc}$ has the ability to replicate by splitting, we
assume binary splitting here.

So there are three main processes which govern the dynamics of
prions in this model.
\begin{itemize}
\item growth in length by polymerization with rate $\tau>0$;
\item binary splitting with rate $\beta(x)>0$, a polymer of length $x>0$ splits into one of length $0<y<x$ and
one of length $x-y$ with probability $\kappa(y,x)$;
\item natural degradation with rate $\gamma>0$ for the monomers and with rate $\mu(x)$ for the polymers with length $x$.
\end{itemize}

The model proposed in \cite{NKKM98} further assumes that polymers
of length $0<x\leq x_0$ immediately decompose completely into
monomers. This reflects the assumption that $PrP^{Sc}$ polymers
are unbranched and form a simple $\alpha$-helix with $x_0$ monomer
units per turn. An $\alpha$-helix of length less than $x_0$ is
incomplete and thus is much less stable. Denoting the numbers of
monomers at time $t$ by $V(t)$ and the density of polymers by
$u(t,x)$, we obtain the following model equations.
\begin{eqnarray}
\label{pde}
&&\partial_tV(t) =\lambda -\gamma V(t) -\tau V(t)\int_{x_0}^\infty u(t,x)dx +2\int_0^{x_0}x\int_{x_0}^\infty\beta(y)\kappa(x,y)u(t,y)dydx\nonumber\\
&&\partial_t u(t,x)+\tau V(t)\partial_x u(t,x)+(\mu(x)+\beta(x))u(t,x)=2\int_x^\infty\beta(y)\kappa(x,y)u(t,y)dy\\
&& V(0)=V_0\geq 0, \quad u(t,x_0)=0, \quad u(0,x)=u_0(x),\nonumber
\end{eqnarray}
where $t\geq0$ and $x_0\leq x <\infty$. Here $\lambda>0$ is a
constant background source of monomers. Observe that the splitting
function $\kappa(y,x)$ should satisfy the following properties.
$$ \kappa(y,x)\geq0,\quad \kappa(y,x)=\kappa(x-y,x),\quad \int_0^x \kappa(y,x)dy =1,$$
for all $x\geq x_0$, $y\geq0$, and $\kappa(y,x)=0$ if $y>x$ or
$x\leq x_0$. Note that these conditions imply
$$ 2\int_0^x y\kappa(y,x)dy =x,\quad x>0.$$
In fact,
\begin{eqnarray*}
&& 2\int_0^x y\kappa(y,x)dy = \int_0^x y\kappa(y,x)dy + \int_0^x y\kappa(x-y,x)dy\\
&& =\int_0^x y\kappa(y,x)dy + \int_0^x (x-y)\kappa(y,x)dy=
x\int_0^x \kappa(y,x)dy=x.
\end{eqnarray*}

This implies that mass does not change via the splitting process,
and by a simple computation we obtain the following relation for
the total number of monomers in the system.
$$ \frac{d}{dt}[V(t)+ \int_{x_0}^\infty xu(t,x)dx] = \lambda-\gamma V(t) -\int_{x_0}^\infty x\mu(x) u(t,x)dx,\quad t\geq0.$$
In \cite{NKKM98} it is further assumed that  splitting is
equi-distributed (polymer chains are equally likely to split at
all locations), and that the rate of splitting is proportional to
length. This reflects again the hypothesis that polymers form
$\alpha$-helices and are not folded in more complicated
configurations, which would make certain segments of the chain
less likely to split than others. Therefore, we make the further
assumptions
$$ \kappa(y,x) = 1/x\; \mbox{ if } x>x_0 \; \mbox{ and } 0<y<x,\quad \kappa(y,x)=0 \; \mbox{ elsewhere },$$
$\beta(x)=\beta x$ is linear, and $\mu(x)\equiv \mu$ constant.
Then the model contains only 6 parameters, and can even be reduced
to a system of 3 ordinary differential equations. In fact,
introduce the new functions
$$ U(t)=\int_{x_0}^\infty u(t,y)dy \quad \mbox{ and } \quad P(t)=\int_{x_0}^\infty yu(t,y)dy,$$
representing the total number of polymers, and the total number of
monomers in polymers at time $t$, respectively. Integrating the
equation for $u(t,x)$ over $[x_0,\infty)$ we  get
\begin{eqnarray*}
\frac{d}{dt} U(t)&=&-\tau V(t)u(t,x)|_{x_0}^\infty-\mu U(t)-\beta P(t)+2\beta \int_{x_0}^\infty\int_x^\infty u(t,y)dydx\\
&=& -\mu U(t)-\beta P(t) +2\beta\int_{x_0}^\infty u(t,y)(y-x_0)dy\\
&=& -\mu U(t)-\beta P(t)+2\beta P(t)-2\beta x_0 U(t),
\end{eqnarray*}
hence
$$\dot{U}(t)= -(\mu+2\beta x_0) U(t)+\beta P(t).$$
Multiplying the equation for $u(t,x)$  by $x$, integration yields
\begin{eqnarray*}
\frac{d}{dt} P(t)&=&-\tau V(t)(xu(t,x)|_{x_0}^\infty-\int_{x_0}^\infty u(t,y)dy)\\
&&-\mu P(t)-\beta \int_{x_0}^\infty u(t,x)x^2dx+2\beta \int_{x_0}^\infty x\int_x^\infty u(t,y)dydx\\
&=& \tau V(t)U(t)-\mu P(t)-\beta \int_{x_0}^\infty u(t,x) x^2dx +\beta\int_{x_0}^\infty u(t,y)(y^2-x_0^2)dy\\
&=& \tau V(t)U(t) -\mu P(t) -\beta x_0^2 U(t),
\end{eqnarray*}
hence
$$\dot{P}(t)= \tau U(t)V(t)-\mu P(t)-\beta x_0^2U(t).$$

Thus we obtain the following closed model involving only ordinary
differential equations.
\begin{eqnarray}
\label{model} \dot{U}&=& \beta P -\mu U -2\beta x_0  U\nonumber\\
\dot{V} &=& \lambda -\gamma V -\tau UV +\beta x_0^2 U\\
\dot{P} &=& \tau UV -\mu P -\beta x_0^2 U\nonumber
\end{eqnarray}
with initial conditions
$$ U(0)=U_0\geq 0,\quad V(0)=V_0\geq 0, \quad P(0)=P_0\geq x_0 U_0.$$
This way the partial differential equation for the density
$u(t,x)$ decouples from the ordinary differential equations. Once
the solutions of \eqref{model} are known, one has to solve only a
linear partial integro-differential eqution to obtain $u(t,x)$.
The system (1.2) is identical to the "basic virus dynamics model"
that is discussed at length in \cite{MaNo02}.

Concerning the ode-system \eqref{model} we have the following
result from Pr\"uss, Pujo-Menjouet, Webb and Zacher \cite{PPWZ05}.

\begin{theorem} Suppose $x_0,\beta,\gamma,\lambda,\mu,\tau>0$ are given constants. Then the system
(\ref{model}) induces a global semiflow on the set
$K=\{(U,V,P)\in\R^3:\; U,V,P-x_0U\geq0\}$. There is precisely one
disease free equilibrium $(0,\lambda/\gamma,0)$ which is globally
exponentially stable if and only if
$\mu+x_0\beta>\sqrt{\lambda\beta\tau/\gamma}$, and asymptotically
stable in case of equality. On the other hand, if
$\mu+x_0\beta<\sqrt{\lambda\beta\tau/\gamma}$ there is the unique
disease equilibrium
$$\Big(\frac{\lambda\beta\tau-\gamma(\mu+\beta
x_0)^2}{\mu\tau(\mu+2\beta x_0)},\frac{(\mu+\beta
x_0))^2}{\beta\tau}, \frac{\lambda\beta\tau-\gamma(\mu+\beta
x_0)^2}{\beta\mu\tau}\Big)$$
which is globally exponentially
stable in $K\setminus[\{0\}\times\R_+\times\{0\}]$.
\end{theorem}

It is the purpose of this paper to study the full system
\eqref{pde} under the assumptions of equi-distributed splitting,
linear splitting rate, and constant rates of degradation.

Since $V(t)+\int_{x_0}^\infty xu(t,x)dx$ is the total number of
monomers in the system, which should be finite at any time, it
seems reasonable to study \eqref{pde} in the standard cone $Z_+:=
\R_+\times L_1^+((x_0,\infty);xdx)$ of the Banach space $Z:=
\R\times L_1((x_0,\infty);xdx)$. The following theorem summarizes
our results.

\begin{theorem}
Assume equi-distributed splitting with linear splitting rate
$\beta(x)=\beta x$ and constant degradation rates $\gamma$ and
$\mu(x)\equiv \mu$. Suppose $\lambda,\tau,\beta,\gamma,\mu,x_0>0$.
Then \eqref{pde} generates a global semiflow in the natural phase
space $Z_+$. Furthermore,
\\(i) \, if $\lambda\beta\tau/\gamma\leq (\mu+\beta x_0)^2$, then
the disease-free equilibrium $\bar{z}=(\lambda/\gamma,0)$ is
globally  asymptotically stable in $Z_+$,  and even exponentially
in the case  of strict inequality;
\\(ii) \, if $\lambda\beta\tau/\gamma> (\mu+\beta x_0)^2$, then there is a
unique disease equilibrium $z_*=(V_*,u_*)$ which is globally
asymptotically stable in $Z_+\setminus(\R_+\times\{0\})$. It is
given by
$$ V_*= \frac{(\mu+\beta x_0)^2}{\beta\tau},\quad u_*(x)= \frac{2\beta}{\mu\tau}
\frac{\lambda\beta\tau-\gamma(\mu+\beta x_0)^2}{(\mu+\beta
x_0)(\mu+2\beta x_0)}\Phi\big(\frac{\beta(x-x_0)}{\mu+\beta
x_0}\big),$$
where $\Phi(r)= (r+r^2/2)\exp(-(r+r^2/2))$.
\end{theorem}

The remaining part of this paper deals with the proof of this
result. Recall that the function $\omega(t):= \tau V(t)$ can be
considered as known, by Theorem 1.1, and $\omega(t)\rightarrow
\omega_\infty$ exponentially, where either $\omega_\infty =
\lambda/\gamma$ in the disease-free or $\omega_\infty= (\mu+\beta
x_0)^2/\beta$ in the disease case. Hence we have to solve a linear
nonautonomous partial integro-differential equation of first
order. For this we shall use standard techniques from the theory
of $C_0$-semigroups and we refer to the monograph Arendt, Batty,
Hieber and Neubrander \cite{ABHN01} as a general reference for the
results employed below.

We proceed in four  steps. First we study the autonomous case
where $\omega\equiv \omega_\infty$. In Section 2 we show that
there is a unique $C_0$-semigroup $T(t)=e^{-Lt}$ associated with
the pde-part of \eqref{pde} in $X=L_1((x_0,\infty);xdx)$, which is
positive and contractive, and even exponentially stable in the
disease-free case. The resolvent of $L$ is shown to be compact in
Section 3, hence $L$ has only point spectrum in the closed right
half-plane. In the disease case, we further show that 0 is the
only eigenvalue of $L$ on the imaginary axis, it is simple and so
the ergodic projection $\cP$ onto the kernel $N(L)$ of $L$ along
the range $R(L)$ of $L$ exists and is rank one. We compute an
element $e\in N(L)$ which is positive. A result of Arendt, Batty,
Lubich and Phong \cite{ABHN01} then shows that $T(t)$ is strongly
ergodic, i.e. $\lim_{t\rightarrow\infty}T(t)=\cP$ strongly in $X$.
Wellposedness of the nonautonomous problem is proved in Section 4
by means of  monotone convergence, it is shown that the evolution
operator exists and is bounded. Moreover, bounds for $\partial_x
u(t,\cdot)$ in $X$ are derived. Finally, in Section 5 we put
together these results to prove Theorem 1.2.

While we assume throughout that $\beta(x) = \beta x, \, \mu(x) =
\mu$ (constant), and $y\kappa(x,y) = 1$ for $x < y ,\, y > x_0$,
$\kappa(x,y) = 0$ elsewhere, our methods extend to versions of
(1.1) where these assumptions do not hold. We do not carry out
these generalizations since it is not clear which would be
biologically reasonable. On the other hand, the equation discussed
in this paper
$$
\partial_t u(t,x) = -\tau V(t)\partial_x u(t,x) -(\mu +\beta x)u(t,x) +2\beta\int_x^\infty u(t,y)dy
$$
for $x > x_0, \, t> 0$, with initial and boundary data as in
(1.1), can be solved with an integral transformation followed by
the method of characteristics. Namely, define
$$v(t,x)= \int_x^\infty \int_y^\infty u(t,\xi) \, d\xi \, dy =
\int_x^\infty (\xi - x)u(t,\xi) \, d\xi, \quad \partial_x^2 v(t,x)
= u(t,x) \, .
$$
Then a computation shows that $v$ solves the first order partial
differential equation without integral term
$$
\partial_t v(t,x)=-\tau V(t) \partial_x v(t,x)-(\mu+\beta x)v(t,x)
$$
for $x > x_0, \, t>0$, with initial data $v(0,x)$ obtained by
integrating $u_0$ twice and boundary data $v(t,x_0) =  P(t) - x_0
U(t)$. The equation for $v$ may be solved by the method of
characteristics, and $u$ is recovered from $\partial_x^2 v(t,x) =
u(t,x)$. The solution depends on the initial data in the region
$\{(x,t) \, | \, x > x_0 + \tau \int_0^tV(s) ds \, \}$ and on the
boundary data in the complement of this region. Since $V(t)$
always has a positive limit, it is evident that the contribution
from the initial data is swept out towards large $x$-values and
decays exponentially, in fact, at a rate like $e^{-\epsilon t^2}$
for some $\epsilon > 0$. If the disease-free state is stable, then
$\left( P(t), U(t) \right) \to (0,0)$ as $t \to \infty$, which
implies that the solution $u$ converges to zero also in the region
where it depends on the boundary data. In the case of a positive
disease equilibrium, $P(t) - x_0 U(t)$ has a positive limit as $t
\to \infty$, which determine the limiting equilibrium distribution
$u_*$ given in Theorem 1.2. This method breaks down if
$\beta(\cdot), \, \mu(\cdot)$, or $\kappa(\cdot,\cdot)$ have more
complicated forms, as the reader will readily confirm.

\section{The Linear Automomous Problem}

\subsection{Functional Analytic Setting}

We consider the problem
\begin{eqnarray}
\label{pde0}
\partial_t u(t,x)+\omega \partial_x u(t,x)+(\mu+\beta x)u(t,x)= 2\beta\int_x^\infty u(t,y)dy,\\
u(0,x)= u_0(x),\quad u(t,x_0)=0,      \quad t>0,\;x>x_0.\nonumber
\end{eqnarray}
Set $w(t,x)=u(t,x+x_0)$, $x\geq0$. Then this problem becomes the
following one on $\R_+$.
\begin{eqnarray}
\label{pde1}
\partial_t w(t,x)+\omega \partial_x w(t,x)+(\mu_0+\beta x)w(t,x)= 2\beta\int_x^\infty w(t,y)dy,\\
w(0,x)=g(x):=u_0(x+x_0),\quad w(t,0)=0,      \quad
t>0,\;x>0.\nonumber
\end{eqnarray}
Here we have set $\mu_0=\mu+\beta x_0$. $\omega$ plays the role of
$\tau V$ at $\infty$, i.e.\
$$\omega=\tau V(\infty)=\lambda\tau/\gamma$$ in the disease-free case or
$$\omega=\tau V(\infty)=(\mu+\beta x_0)^2/\beta=\mu_0^2/\beta$$
in the disease case.

We want to study (\ref{pde1}) in the basic space
$X=L_1(\R_+;(a+x)dx)$, where we choose as the norm
$$ ||w||=a|w|_1+|xw|_1,$$
with $a>0$ to be determined later. We define two linear operators
in $X$ by means of
$$ Au(x)= \omega u^\prime(x)+(\mu_0+\beta x)u(x),\quad x\in\R_+,$$
with domain
$$D(A)=\{ u\in W^1_1(\R_+)\cap X: \; x^2u\in L_1(\R_+), x u^\prime(x)\in L_1(\R_+),\; u(0)=0\},$$
and $$Bu(x)=2\beta\int_x^\infty u(y)dy,\quad D(B)=D(A).$$ Both
operators are well-defined and linear, $B$ will be considered as a
perturbation of $A$.

\subsection{$m$-Accretivity of $A$}

We have
\begin{eqnarray*}
\int_0^\infty Au \sgn{u} dx &=& \omega\int_0^\infty |u|^\prime dx +\mu_0|u|_1 +\beta|xu|_1\\
&=& \mu_0|u|_1 +\beta|xu|_1,
\end{eqnarray*}
and
\begin{eqnarray*}
\int_0^\infty Au \sgn{u} xdx &=& \omega\int_0^\infty |u|^\prime xdx +\mu_0|xu|_1 +\beta|x^2u|_1\\
&=& -\omega|u|_1+\mu_0|xu|_1 +\beta|x^2u|_1.
\end{eqnarray*}
Employing the bracket in $L_1$ this implies
$$[Au,u]_+\geq (a\mu_0-\omega)|u|_1+(a\beta +\mu_0)|xu|_1\geq \eta||u||,$$
for some $\eta>0$ provided $\mu_0>\omega/a$. Hence for such $a$,
$A$ is strictly accretive, in particular closable.

Next we compute the resolvent of $A$. The equation
$(\lambda+A)u=f$ is equivalent to solving the ode
\begin{equation}\label{help}
\lambda u(x)+ \omega u^\prime(x)+(\mu_0+\beta x)u(x)=f(x),\quad
x>0,
\end{equation}
with initial condition $u(0)=0$. Therefore we obtain
$$u=(\lambda+A)^{-1} f(x)=\frac{1}{\omega}\int_0^x
\exp{-[(\lambda+\mu_0)(x-y)/\omega+\beta(x^2-y^2)/2\omega]}
f(y)dy.$$

If $f\in L_1(\R_+)$ then on easily obtains the estimate
$$ |u|_1\leq |f|_1/(\lambda+\mu_0).$$
If also $xf\in L_1(\R_+)$ then
\begin{eqnarray*}
|x^2u(x)|&\leq& \frac{1}{\omega}\int_0^x e^{-(\lambda+\mu_0)(x-y)/\omega}(x^2-y^2)e^{-\beta(x^2-y^2)/2\omega)}|f(y)|dy\\
&+& \frac{1}{\omega}\int_0^x ye^{-\beta(x-y)2y/2\omega}y|f(y)|dy,
\end{eqnarray*}
hence
$$|x^2u|_1\leq \frac{1}{\omega}\frac{\omega}{\lambda+\mu_0}\frac{2\omega}{\beta e}|f|_1
+\frac{1}{\omega}\frac{\omega}{\beta^2}|xf|_1.$$
This shows that
$x^2u\in L_1(\R_+)$, hence  $xu\in L_1(\R_+)$, and then by
equation \eqref{help} also $u^\prime\in L_1(\R_+)$ as well as $x
u^\prime\in L_1(\R_+)$, i.e.\ $u\in D(A)$. This shows that $A$ is
$m$-accretive.

As a consequence we note that $-A$ generates a $C_0$-semigroup in
$X$ which is also positive and strictly contractive, hence
exponentially stable.

\subsection{Accretivity of $A-B$}

We have
$$|\int_x^\infty u(x)dx|_1\leq |xu|_1,\quad |x\int_x^\infty u(x)dx|_1\leq \frac{1}{2}|x^2u|_1,$$
and therefore
$$\int_0^\infty (Au-Bu)\sgn(u)dx\geq \mu_0|u|_1+\beta|xu|_1-2\beta|xu|_1,$$
as well as
$$\int_0^\infty(Au-Bu)\sgn(u) xdx\geq -\omega|u|_1+\mu_0|xu|_1.$$
This yields
$$[(A-B)u,u]_+\geq (\mu_0a-\omega)|u|_1+(\mu_0-\beta a)|xu|_1\geq0,$$
for all $u\in D(A)$, provided $\mu_0a\geq\omega$ and $\mu_0\geq
\beta a$. Such a choice of $a>0$ is possible if and only if the
condition $\omega/\mu_0\leq \mu_0/\beta$ is met, i.e.\ if and only
if
$$\omega\leq \mu_0^2/\beta$$
holds true. Now in the disease-free case  we have
$\omega=\lambda\tau/\gamma$, while in the disease case
$\omega=\mu_0^2/\beta$; then $a=\mu_0/\beta$. Thus $A-B$ will be
strictly accretive in the disease-free case while it will be
accretive only in the disease case. In the first case, the decay
rate can easily be estimated not to be smaller than
$\mu_0-\sqrt{\lambda\beta\tau/\gamma}$.

\subsection{Density of the Range of $A-B$}

Let $f\in L_1(\R_+;(a+x)dx)$ be given and assume $f\geq0$. Set
$u_1=(1+A)^{-1}f$ and define the sequence $u_n$ inductively by
means of
$$ u_{n+1}=u_1+(1+A)^{-1} Bu_n.$$
Then $u_1\geq0$, and $u_{2}-u_1=(1+A)^{-1}Bu_1\geq0$, hence by
induction $u_{n+1}\geq u_n$ pointwise, since $B$ is positive. This
shows that the sequence of functions $u_n$ is nonnegative and
increasing pointwise. Moreover,
$$\omega u^\prime_n +(1+\mu_0+\beta x)u_n= f+2\beta\int_x^\infty u_{n-1}(y)dy\leq f+2\beta\int_x^\infty u_n(y)dy,$$
which implies
$$(1+\mu_0)|u_n|_1+\beta|xu_n|_1\leq |f|_1+2\beta|xu_n|_1,$$
and
$$ -\omega|u_n|_1+(1+\mu_0)|xu_n|_1+\beta^2|x^2u_n|_1\leq |xf|_1+\beta|x^2u_n|_1.$$
Choosing $a$ as above this yields an a priori bound for the
sequence $(u_n)$
$$||u_n||=a|u_n|_1+|xu_n|\leq C||f||,$$
and  therefore we may conclude by the monotone convergence theorem
$u_n\rightarrow u_\infty$ as $n\rightarrow\infty$. If in addition
$x^2f\in L_1(\R_+)$ then we obtain in a similar way boundedness of
$x^2u_n$ in $X$. This implies $(1+A-B)u_n=
f+B(u_{n-1}-u_n)\rightarrow f$ in $X$ as $n\rightarrow\infty$,
hence $u_\infty\in D(\overline{A-B})$ and
$u_\infty=(1+\overline{A-B})^{-1}f$. Since $L_1=L_1^+-L_1^+$ we
may conclude $R(1+\overline{A-B})=X$, i.e.\ the closure of $A-B$
is $m$-accretive.

\bigskip

\begin{remark}
The above proof shows that the resolvent of $\overline{A-B}$ is
positive, hence the semigroup generated by this operator will be
as as well.
\end{remark}

\subsection{ Irreducibility}

Suppose $f\in X$ is nonnegative and $u$ solves
$$\omega u^\prime +(\lambda +\mu_0+\beta x)u =f+2\beta\int_x^\infty u(y)dy,\quad x\geq0,$$
with initial value $u(0)=0$. If $f\not\equiv0$ then let
$x_1:=\inf\supp f$. We have
$$u(x)=\frac{1}{\omega}\int_0^x \exp{-[(\lambda+\mu_0)(x-y)/\omega+\beta(x^2-y^2)/2\omega]} [f(y)+Bu(y)]dy.$$
Since we already know $u(x)\geq0$, this formula implies $u(x)>0$
for all $x>x_1$. But then $\int_x^\infty u(y)dy>0$ for all
$x\geq0$, and so so $u(x)>0$ for all $x>0$. This proves the
irreducibility of the semigroup generated by $\overline{A-B}$.

\subsection{$A-B$ is not Closed}

Unfortunately, the sum $A-B $ is not closed. We show this by the
following example.

\begin{example}
Set $u=\chi/x^3$ where $\chi$ denotes a cut-off function which is
$0$ on $[0,1]$ and $1$ on  $[2,\infty)$. Then $u,u^\prime u,xu\in
L_1(\R_+)$, but $x^2u\not\in L_1(\R_+)$, and $u(0)=0$. On the
other hand,
\begin{eqnarray*}
f(x)&:=& \omega u^\prime(x)+(\lambda +\mu_0+\beta x)u(x)-2\beta\int_x^\infty u(y)dy\\
&=& \omega\chi^\prime/x^3-3\omega\chi/x^4+(\lambda+\mu_0)\chi/x^3
+\beta \chi/x^2- 2\beta \int_x^\infty\chi(y)dy/y^3
\end{eqnarray*}
Since
\begin{eqnarray*}
&&\chi(x)/x^2-2\int_x^\infty\chi(y)dy/y^3 = \chi(x)/x^2+\chi(y)/y^2|_x^\infty-\int^\infty_x\chi^\prime(y)dy/y^2\\
&=&-\int_x^\infty \chi^\prime(y)dy/y^2,
\end{eqnarray*}
we obtain
$$f=\omega \chi^\prime(x)/x^3-3\omega\chi(x)/x^4+(\lambda+\mu_0)\chi(x)/x^3-\beta \int_x^\infty \chi^\prime(y)dy/y^2.$$
Obviously, $f$ as well as $xf$ belong to $L_1(\R_+)$, so $A-B$
with domain $D(A)$ is not closed.
\end{example}

\subsection{Summary}

Let us summarize what we have shown so far.

\begin{theorem}
Suppose $\beta\omega\leq \mu_0^2$. Then problem (\ref{pde1}) is
well-posed in $X=L_1(\R_+;(a+x)dx)$ and admits an associated
$C_0$-semigroup $T(t)=e^{-Lt}$ which is positive. If $a$ is chosen
from the interval $a\in[\omega/\mu_0, \mu_0/\beta]$ then $T(t)$ is
nonexpansive.

In the strictly disease free case
$\omega=\lambda\tau/\gamma<\mu_0^2/\beta$, the semigroup $T(t)$ is
exponentially stable with type $\omega_0(T)\leq
-\mu_0+\sqrt{\lambda\beta\tau/\gamma}<0$.
\end{theorem}

\section{Asymptotic Behavior of the Autonomous Problem}

\subsection{Compactness}

Set $L=\overline{A-B}$. Since $L$ is m-accretive in
$X=L_1(\R_+;(a+x)dx)$, the spectrum $\sigma(L)$ is contained in
the closed right halfplane. We want to show that the resolvent of
$L$ is compact. For this purpose we derive another representation
of $(\lambda+L)^{-1}$ for $\lambda>0$. Let $f\in X$ and set
$u=(\lambda+L)^{-1}f$. Then we obtain
$$u=(\lambda+A)^{-1}f + (\lambda+A)^{-1} B u,$$
and
\begin{eqnarray*}
(\lambda+A)^{-1}Bu&=& 2\beta(\lambda+A)^{-1}[\int_x^\infty u(y)dy]\\
&=& \frac{2\beta}{\omega} \int_0^x e^{-(\lambda +\mu_0)(x-y)/\omega}e^{-\beta(x^2-y^2)/2\omega}[\int_y^\infty u(r)dr]dy\\
&=& \frac{2\beta}{\omega} \int_x^\infty u(r)[\int_0^x e^{-(\lambda +\mu_0)(x-y)/\omega}e^{-\beta(x^2-y^2)/2\omega}dy]dr\\
&+& \frac{2\beta}{\omega} \int_0^x u(r) [\int_0^r e^{-(\lambda +\mu_0)(x-y)/\omega}e^{-\beta(x^2-y^2)/2\omega}dy]dr\\
&=& k_\lambda(x) \int_x^\infty u(r)dr + \omega(\lambda+A)^{-1}
[k_\lambda u],
\end{eqnarray*}
where
$$k_\lambda(x)= \frac{2\beta}{\omega}\int_0^x e^{-(\lambda +\mu_0)(x-y)/\omega}e^{-\beta(x^2-y^2)/2\omega}dy.$$
Note that
$$0\leq k_\lambda(x)\leq \frac{2\beta}{\omega}\int_0^x e^{-(\lambda +\mu_0)(x-y)/\omega}dy\leq \frac{2\beta}{\lambda+\mu_0},$$
i.e.\ $k_\lambda\in L_\infty(\R_+)$. We thus have the identity
$$u(x)-k_\lambda(x)\int_x^\infty u(y)dy = (\lambda+A)^{-1}f(x)+ \omega(\lambda+A)^{-1}[k_\lambda u]=: g(x),$$
and $u(0)=0$. We may solve this equation for $u$ to the result
$$ u(x)= g(x)-k_\lambda(x)\int_0^x \exp\left(-\int_y^x k_\lambda(r)dr\right) g(y)dy + k_\lambda(x) \exp \left(-\int_0^x k_\lambda(s)ds\right)<q_\lambda| f>,$$
where
$$<q_\lambda,f> :=
\frac{1}{(\lambda+\mu_0)^2-\omega\beta}((\lambda+\mu_0)\int_0^\infty
f(s)ds+\beta \int_0^\infty sf(s)ds)\,.$$
This way we have the
representation
\begin{equation}
\label{repres} (\lambda+L)^{-1}f = (1-R_\lambda)(\lambda+A)^{-1}[
1+\omega k_\lambda (\lambda+L)^{-1}]f +k_\lambda(x)
\exp\left(-\int_0^x k_\lambda(s)ds \right)<q_\lambda| f>,
\end{equation}
with
$$(R_\lambda g)(x) = k_\lambda(x)\int_0^x \exp \left(-\int_y^x k_\lambda(r)dr\right) g(y)dy.$$
Next $D(A)$ embeds compactly into $X$, hence $(\lambda+A)^{-1}$ is
compact. From boundedness of $k_\lambda$ we may then conclude that
$(\lambda+L)^{-1}$ is compact, as soon as we know that the
Volterra operator $R_\lambda$ is bounded in $X$.

To prove the latter we estimate as follows
\begin{eqnarray*}
 ||R_\lambda g||&=&\int_0^\infty (a+x) k_\lambda(x)|\int_0^x \exp\left(-\int_y^x k_\lambda(r)dr\right) g(y)dy|dx\\
&\leq& \int_0^\infty|g(y)|[\int_y^\infty (a+x)k_\lambda(x)\exp\left(-\int_y^x k_\lambda(r)dr\right)dx]dy\\
&=& \int_0^\infty |g(y)| [ (a+y) + \int_y^\infty \exp\left(-\int_y^x k_\lambda(r)dr\right)dx]dy\\
&\leq& C_\lambda\int_0^\infty |g(y)|(a+y)dy=C_\lambda||g||,
\end{eqnarray*}
as we show now.
\begin{eqnarray*}
k_\lambda(x)&=& \frac{2\beta}{\omega}\int_0^x e^{-(\lambda +\mu_0)(x-y)/\omega}e^{-\beta(x^2-y^2)/2\omega}dy\\
&\geq & \frac{2\beta}{\omega}\int_0^x e^{-(\lambda +\mu_0)y/\omega}e^{-\beta xy/\omega}dy\\
&=& \frac{2\beta}{\lambda+\mu_0+\beta x} (1- e^{-(\lambda+\mu_0+\beta x)x/\omega})\\
&\geq & \frac{2\beta}{\lambda+\mu_0+\beta x}\cdot \frac{(\lambda+\mu_0+\beta x)x/\omega}{1+(\lambda+\mu_0+\beta x)x/\omega)}\\
&=& \frac{2\beta x}{\omega +(\lambda+\mu+\beta x)x},
\end{eqnarray*}
by the elementary inequality $1-e^{-x} \geq x/(1+x)$. This implies
\begin{eqnarray*}
\int_y^x k_\lambda(r)dr &\geq& 2\beta\int_y^x rdr/(\omega +(\lambda+\mu_0+\beta r)r\\
&=& \int_y^x \frac{2\beta r
+\lambda+\mu_0}{\omega+(\lambda+\mu_0)r +\beta r^2}dr
-(\lambda+\mu_0)\int_y^x \frac{dr}
{\omega+(\lambda+\mu_0)r +\beta r^2}\\
&\geq& \log \frac{\omega +(\lambda+\mu_0)x + \beta
x^2}{\omega+(\lambda+\mu_0)y +\beta y^2} - c_\lambda,
\end{eqnarray*}
since the second integral is bounded. This estimate finally yields
$$\int_y^\infty \exp\left(-\int_y^x k_\lambda(r)dr\right)dx\leq e^{c_\lambda}\int_y^\infty\frac{\omega +(\lambda+\mu_0)y + \beta y^2}
{\omega+(\lambda+\mu_0)x +\beta x^2}dx\leq C_\lambda(a+y).
$$
This completes the proof of compactness of the resolvent of $L$.

\subsection{Ergodicity}

Since the resolvent of $L$ is compact we know that the spectrum of
$L$ consists only of eigenvalues of finite multiplicity, these are
poles of the resolvent of $L$. By accretivity of $L$ we have the
inequality $|(\lambda+L)^{-1}|_{{\cal B}(X)}\leq 1/{\rm
Re}\lambda$, ${\rm Re} \lambda>0$, hence the resolvent can only
have poles of first order on the imaginary axis. This shows that
all eigenvalues on the imaginary axis are semisimple. Compactness
of the resolvent implies also that the range of $\lambda+L$ is
closed, for each $\lambda\in\C$. In particular, we have the direct
sum decomposition $X= N(L)\oplus R(L)$, i.e.\ ergodicity in the
sense of Abel.

Now we concentrate on the disease equilibrium which means
$a=\mu_0/\beta$ and $\omega=\mu_0^2/\beta$. A function $e(x)$
belongs to the kernel of $L$ if
$$ \omega e^\prime (x)+(\mu_0+\beta x) e(x)- 2\beta \int_x^\infty e(y)dy =0, \quad x>0,\; e(0)=0,$$
or equivalently
$$ e^{\prime\prime}(x)+\frac{\beta}{\mu_0}(1+\frac{\beta}{\mu_0}x) e^\prime(x)+ 3 \frac{\beta^2}{\mu_0^2} e(x)=0, \quad x>0,\; e(0)=0.$$
The scaling $e(x)= v(\beta x/\mu_0)$ reduces this problem to
$$ v^{\prime\prime}(z) +(1+z)v^\prime(z)+3 v(z)=0,\quad z>0,\; v(0)=0.$$
By the initial condition $v(0)=0$, this shows that the kernel of
$L$ can be only one-dimensional, and a simple computation yields
that
$$v(z)= (z+z^2/2) e^{-(z+z^2/2)}, \quad z>0,$$
is a solution. Therefore $N(L)={\rm span}\{ e\}$, with $e(x)=
(\beta/\mu_0)^2v(\beta x/\mu_0)$, and another simple computation
yields
$$ \int_0^\infty (a+x) e(x)dx =1.$$
Since $L$ is Fredholm with index zero, the kernel $N(L^*)$ of the
dual of $L$ has also a one-dimensional kernel which are the
constant functions. The ergodic projection $\cP$ onto the kernel
of $L$ along the range of $L$ is then given by
\begin{equation}\label{projection}
\cP u(x)= [\int_0^\infty(a+x) u(x)dx] e(x)= <u|e^*>e(x), \quad
x>0.
\end{equation}
Suppose there are no other eigenvalues of $L$ on the imaginary
axis. Then $L^*$ also has no other eigenvalues on the imaginary
axis, and then by the theorem of Arendt, Batty, Lubich and Phong
we may conclude that
$$ e^{-Lt} u \rightarrow \cP u \quad \mbox{ as } t\rightarrow\infty, \mbox{ for each } u\in X,$$
i.e.\ the semigroup generated by $-L$ is strongly ergodic.

We show now that there are in fact no eigenvalues other than 0 on
the imaginary axis. Suppose on the contrary that
$$ i\rho u(x) +\omega u^\prime(x) +(\mu_0+\beta x)u(x)= 2\beta\int_x^\infty u(y)dy, \quad x>0,\; u(0)=0,$$
$u\neq0$. Multiplying this equation with $\bar{u}/|u|$, taking
real parts, and integrating over $\R_+$ we obtain
\begin{equation}\label{i}
\mu_0|u|_1 +\beta|xu|_1 = 2\beta {\rm Re}\int_0^\infty
u(x)\int_0^x \bar{u}(y)/|u(y)|dydx\leq 2\beta |xu|_1,
\end{equation}
and similarly, multiplying with $x\bar{u}(x)/|u(x)|$ we get
\begin{equation}\label{ii}
-\omega|u|_1+ \mu_0|xu|_1+\beta |x^2u|_1= 2\beta {\rm Re}
\int_0^\infty u(x)\int_0^x y\bar{u}(y)/|u(y)|dydx\leq \beta
|x^2u|_1.
\end{equation}
Multiplying the first inequality with $a=\mu_0/\beta$ and adding
the second we arrive at a contradiction if at least one of the
inequalities (\ref{i}), (\ref{ii}) is strict. Hence we must have
$$ {\rm Re}\int_0^\infty u(x)\int_0^x \bar{u}(y)/|u(y)|dydx= |xu|_1,$$
which implies with $\arg u(x)=\theta(x)$
$$x\equiv{\rm Re}\int_0^x e^{i(\theta(x)-\theta(y))}dy= \frac{1}{2}\frac{d}{dx} |\int_0^x e^{i\theta(y)}dy|^2,$$
or equivalently
$$ |\int_0^x e^{i\theta(y)}dy|^2=x^2, \quad x>0.$$
But this is only possible if $\theta(y)$ is constant, w.l.o.g.\ we
may assume $\theta =0$ i.e.\ $u(x)$ is nonnegative, which in turn
yields $\rho=0$ since $u\neq0$ by assumption.

\bigskip

\subsection{Summary}

Let us summarize what we have shown in this section.

\begin{theorem}
Assume the disease case $\omega=\mu_0^2/\beta$, $a=\mu_0/\beta$.
The the semigroup $T(t)=e^{-Lt}$ is strongly ergodic, it converges
strongly to the projection $\cP$ onto the kernel $N(L)$ of $L$
along its range $R(L)$. The kernel is one-dimensional and spanned
by $e(x)= (\beta/\mu_0)^2\Phi(\beta x/\mu_0)$, where
$\Phi(z)=(z+z^2/2)e^{-(z+z^2/2)}$, and the projection $\cP$ is
given by
$$ \cP u(x)=[\int_0^\infty(a+y)u(y)dy]e(x)= <e^*|u>e(x), \quad x>0,\; u\in X.$$
\end{theorem}

\noindent

{\em Remark.} We do not know  whether the ergodicity is
exponential since it is not clear that the type of the semigroup
$e^{-Lt}$ restricted to $R(L)$ is negative.

\section{Well-posedness of the Non-Autonomous Evolution}

\subsection{The Trivial Evolution}

Let $\omega\in C(\R_+)$ be positive, such that
$0<\omega_\infty=\lim_{t\rightarrow\infty} \omega(t)$ exists, and
assume $\omega(\cdot)-\omega_\infty\in L_1(\R_+)$. Let
$$\omega_+=\max_{s\geq0} \omega(s)\quad \mbox{ and  }\quad \omega_-=\min_{s\geq0} \omega(s),$$
and note that $\omega_+\geq \omega_->0$. We are particularly
interested in the cases $\omega_\infty=\lambda\tau/\gamma$, the
disease-free case, and $\omega_\infty= \mu_0^2/\beta$, the disease
case. We want to show that the nonautonomous problem is well-posed
in $X=L_1(\R_+;(a+x)dx)$. We begin with the problem
\begin{eqnarray}
\label{trivial}
\partial_t u(t,x)+\omega(t)\partial_x u(t,x) +(\mu_0+\beta x) u(t,x)=0,\quad x>0, t>s\geq0\\
u(s,x)=g(x),\quad u(t,0)=0,\quad t>s\geq0, x>0.\nonumber
\end{eqnarray}
The method of characteristics yields easily the evolution operator
$U_0(t,s)$ for this problem. It is given by
\begin{eqnarray}
\label{evop}
[U_0(t,s)g](x)=u(t,x)= g(x-\int_s^t\omega(\tau)d\tau)e^{-\phi(t,s,x)},\\
\phi(t,s,x)=\mu_0(t-s)+
\beta(t-s)(x-\int_s^t\omega(\tau)d\tau)+\beta\int_s^t(t-\tau)\omega(\tau)d\tau,\nonumber
\end{eqnarray}
if we extend $g$ trivially to $\R$. We obviously have the estimate
$|U_0(t,s)|_{{\cal B}(X)}\leq e^{-\mu_0(t-s)}$, and $u(t,x)$ is a
strong solution in $X$ if the initial function $g$ belongs to $D$
defined by
$$ D:=\{g\in L_1(\R_+):\;x^2g,g^\prime,xg^\prime\in L_1(\R_+), g(0)=0\}.$$
We also need  the solution of
\begin{eqnarray}
\label{trivial1}
\partial_t u(t,x)+\omega(t)\partial_x u(t,x) +(\mu_0+\beta x) u(t,x)=0,\quad x>0, t>s\geq0\nonumber\\
u(s,x)=0,\quad u(t,0)=h(t),\quad t>s\geq0, x>0.
\end{eqnarray}
Again the method of characteristics applies and yields with
$K(t,x) = \int_{\rho(t,x)}^t(r-\rho(t,x))\omega(r)dr$ the formula
$$ [V_0(t,s)h](x)=u(t,x)=h(\rho(t,x))e^{-[\mu_0(t-\rho(t,x)+\beta x(t-\rho(t,x))-\beta K(t,x)]},
$$
for $x<\int_s^t\omega(r)dr$, and zero elsewhere, where the
function $\rho(t,x)$ is defined by the equation
\begin{equation}
\label{kappa} x=\int_\rho^t \omega(r)dr;
\end{equation}
note that this equation has a unique solution $\rho(t,x)\in
(s,t)$, since $\omega(r)\geq \omega_- >0$ for all $r\geq0$, by
assumption, and $x<\int_s^t\omega(r)dr$. Observe that with
$K_0(t,s) =\int_s^t\omega(r)dr$ we have
\begin{eqnarray*}
&&\int_0^\infty (a+x)[V_0(t,s)h](x)dx\leq |h|_\infty\int_0^{K_0(t,s)}(a+x)e^{-\mu_0(t-\rho(t,x))}dx\\
&& \leq |h|_\infty \int_s^t (a+\int_\sigma^t\omega(r)dr)e^{-\mu_0(t-\sigma)} \omega(\rho(t,x)) d\sigma\\
&&\leq |h|_\infty \omega_+\int_0^{t-s}(a+\omega_+\sigma)e^{-\mu_0
\sigma}d\sigma\leq C|h|_\infty ,
\end{eqnarray*}
by the variable transformation $\sigma=\rho(t,x)$. Thus the part
coming from a nontrivial bounded boundary value $h$ is bounded in
$X$.

\subsection{Well-posedness for the Full Problem}

Let us now consider the full problem, i.e.
\begin{eqnarray}
\label{full}
\partial_t u(t,x)+\omega(t)\partial_x u(t,x) +(\mu_0+\beta x) u(t,x)=2\beta\int_x^\infty u(t,y)dy,\\
u(s,x)=g(x),\quad u(t,0)=0,\quad t>s\geq0, \; x>0.\nonumber
\end{eqnarray}
Since the standard cone in $X$ is reproducing, i.e.
$L_1=L_1^+-L_1^+$, we may restrict attention to nonnegative
initial functions $g$. We define the sequence $u_n$ inductively by
$$u_1(t):=U_0(t,s)g, \quad u_{n+1}(t)=u_1(t)+\int_s^t U_0(t,r)Bu_n(r)dr, \quad t\geq s\geq 0.$$
Since $U_0(t,s)$ is positive the functions $u_n$ are as well, and
$u_2(t)\geq u_1(t)$ since $B$ is positive. Inductively we obtain
with
$$u_{n+1}(t)-u_n(t)=\int_s^t U_0(t,r)B(u_n(r)-u_{n-1}(r))dr, \quad t\geq s\geq0,$$
that the functions $u_n$ are pointwise increasing w.r.t.\
$n\in\N$.

Suppose that $g\in D$. Then $u_n$ is a strong solution of
\begin{eqnarray*}
&&\partial_t u_n(t,x)+\omega(t)\partial_x u_n(t,x) +(\mu_0+\beta x) u_n(t,x)=2\beta\int_x^\infty u_{n-1}(t,y)dy\\
&&\qquad\qquad\leq 2\beta\int_x^\infty u_n(t,y)dy,\quad x>0, t>s\geq0\\
&&u(s,x)=g(x),\quad u(t,0)=0,\quad t>s\geq0, x>0,
\end{eqnarray*}
i.e.\ $u_n$ is a strong lower solution of (\ref{full}).
Multiplying the equation with $x^i$ and integrating over $\R_+$
this yields with $z_i(t)=|x^i u_n(t)|_1$
$$\partial_t z_0(t) +\mu_0 z_0(t)+\beta z_1(t)\leq 2\beta z_1(t),$$
for $i=0$, and for $i=1$
$$\partial_tz_1(t)-\omega(t) z_0(t)+\mu_0z_1(t)+\beta z_2(t)\leq \beta z_2(t).$$
Setting $z(t)=(z_0(t),z_1(t))^T$,
$b(t)=(0,(\omega(t)-\omega_\infty)z_0(t))^T$, and defining $G$ by
the $2\times2$-matrix with entries
$-\mu_0,\beta,\omega_\infty,-\mu_0$, this inequality becomes
$$\partial_tz(t)\leq Gz(t)+b(t), \quad t\geq s\geq0.$$
The eigenvalues of $G$ are given by $\lambda_\pm= -\mu_0 \pm
\sqrt{\beta\omega_\infty}$ which are both nonpositive if
$\beta\omega_\infty\leq \mu_0^2$, which is true in both, the
disease-free and the disease case. Since $e^{Gt}$ is positive we
may conclude
$$z(t)\leq e^{G(t-s)}z(s)+\int_s^t e^{G(t-r)}b(r)dr.$$
Boundedness of $e^{Gt}$ then implies an inequality of the form
$$|z(t)|\leq C + C\int_s^t |\omega(r)-\omega_\infty||z(r)|dr,\quad t\geq s\geq0,$$
which implies boundedness of $z(t)$ on $[s,\infty)$ since
$(\omega(\cdot)-\omega_\infty)\in L_1(\R_+)$ by assumption. Note
that the constant $C$ depends only on the parameters
$\mu_0,\beta,\omega_\infty$ and on $||g||$.

Therefore the functions $u_n(t)$ are bounded in $X$ uniformly in
$t$ and $n$. By monotone convergence we may conclude
$u_n(t)\rightarrow u(t)$ in $X$ for each $t\geq s$. Since $B$ is
positive, $Bu_n\rightarrow Bu$ in $L_1(\R_+)$ as well, and then
also
\begin{equation}
\label{mild}
 u(t)=U_0(t,s)g+\int_s^t U_0(t,r)Bu(r)dr,\quad t\geq s\geq0,
\end{equation}
at least in $L_1(\R_+)$. A density argument finally shows that
this conclusion is valid for all initial data $g\in X$.

\bigskip

\noindent

{\em Remark.} It is not clear that solutions of  (\ref{mild}) are
unique. The reason for this  is that $B$ is unbounded. Therefore
we need another definition of mild solution.

\bigskip

\noindent
{\bf Definition.}\,{\em  Let $f\in L_{1,loc}(R_+;X)$. \\
(i)\, We call a  function  $u\in C(\R_+;X)$ strong solution of
\begin{eqnarray}
\label{full1}
\partial_t u(t,x)+\omega(t)\partial_x u(t,x) +(\mu_0+\beta x) u(t,x)=2\beta\int_x^\infty u(t,y)dy+f(t,x),\nonumber\\
u(s,x)=g(x),\quad u(t,0)=0,\quad t>s\geq0, x>0.
\end{eqnarray}
if $u\in C^1(\R_+;X)\cap C(\R_+;D)$ and (\ref{full1}) is valid pointwise.\\
(ii)\, We call a function $u\in C(\R_+;X)$ mild  solution of
(\ref{full1}) if there are $f_n\in L_{1,loc}(R_+;X)$ and strong
solutions $u_n$ of (\ref{full1}) such that $u_n\rightarrow u$ and
$f_n\rightarrow f$ as $n\rightarrow\infty$, in $X$, uniformly on
compact intervals. }

\bigskip

\noindent Suppose that $g\in D$ has compact support. Then each
iteration $u_n(t)$ has also compact support, namely
$$\supp u_n(t)\subset \supp g + \omega_+[0,t],$$
for each $n\in\N$. Therefore each function $u_n(t)$ is a strong
solution of (\ref{full1}) with inhomogeneity
$f_n(t)=B(u_{n-1}(t)-u_n(t))$. This proves that the limit $u(t)$
is a mild solution. Approximation then shows that (\ref{full}) has
at least one mild solution, for each initial value $g\in X$.

Uniqueness of mild solutions can be obtained as follows. If $u$ is
a strong solution of (\ref{full1}) then the equation yields as
above the inequality
$$ \partial_t ||u(t)||\leq \omega_+ ||u(t)||+ ||f(t)||,\quad
t>0,$$
hence
$$ ||u(t)||\leq e^{\omega_+(t-s)}||g||+\int_s^t
e^{\omega_+(t-r)}||f(r)||dr.$$
By approximation this inequality is
also valid for mild solutions, hence $u\equiv 0$ in case $f\equiv
g=0$. Thus mild solutions are unique and of course they satisfy
the integral equation (\ref{mild}).

\subsection{ Summary.}

We have proved the following result about well-posedness of \eqref{full}

\begin{theorem}
Suppose $\omega\in C(\R_+)$ is a given strictly positive function,
such that $\omega_\infty=\lim_{t\rightarrow\infty} \omega(t)>0$
exists and $\omega(\cdot)-\omega_\infty\in L_1(\R_+)$. Then
\eqref{full} is well-posed in the sense of the definition given
above. There exists a unique evolution operator $U(t,s)$ in $X$
generated by \eqref{full}, which is  bounded in $X$, uniformly in
$0\leq s\leq t<\infty$,  and  positive. Moreover, \eqref{full} has
finite speed of propagation with maximum speed less than $
\omega_+= \sup_{t\geq0} \omega(t)$.
\end{theorem}

\subsection{Higher Order Bounds}

Consider an initial function $g\in C_0^\infty(0,\infty)$. Then
$u_1$ is smooth as well and has compact support for each $t\geq
s$. Then the same holds true for $u_2$, hence by induction for all
$u_n$. Setting $v_n=\partial_x u_n$ we have the following problem
for $v_n$.
\begin{eqnarray}
\label{derivative}
&&\partial_t v_n+\omega(t)\partial_x v_n +(\mu_0+\beta x) v_n= -\beta [u_n+2u_{n-1}],\\
&&v_n(s,x)=g^\prime(x),\quad v_n(t,0)=\psi_n(t),\quad t>s\geq0,
x>0\nonumber
\end{eqnarray}
where $\psi_n(t)= \frac{2\beta}{\omega(t)}|u_{n-1}(t)|_1$. This
implies
$$\partial_x u_n(t)=v_n(t) = U_0(t,s)g^\prime -\beta\int_s^t U_0(t,r)[u_n(r)+2u_{n-1}(r)]dr+w_n(t),\quad t\geq s\geq0,$$
with
$$w_n(t)=2\beta V_0(t,s)[|u_{n-1}(\cdot)|_1/\omega(\cdot)].$$
Uniform boundedness of $u_n$ in $X$ and exponential stability of
the evolution operator $U_0(t,s)$ in $X$ then implies boundedness
of $\partial_x u_n$ in $X$. Passing to the limit we get
$$\partial_x u(t) = U_0(t,s)g^\prime -3\beta\int_s^t U_0(t,r)u(r)dr+ w(t),\quad t\geq s\geq0,$$
where
$$w(t,x)=2\beta V_0(t,s)[|u(\cdot)|_1/\omega(\cdot)]. $$
This yields $\partial_xu\in C_b([s,\infty);X)$. The last identity
was proven for $g\in C_0^\infty(0,\infty)$, but via density can be
extended to $g\in D$.

\section{Convergence}

We are now ready to prove the main result on convergence. Let us
first look at the disease-free case. Then with $A(t)$, $B$,
defined as in section 2, and $L(t)=\overline{A(t)-B}$, we know
that $L(t)$ is strictly accretive for large times $t$ if the
parameter $a$ is chosen in
$a\in(\lambda\tau/\gamma\mu_0,\mu_0/\beta)$.  This proves
exponential stability of the trivial solution in the disease-free
case, with decay rate at least
$\mu_0-\sqrt{\lambda\beta\tau/\gamma}$.

Suppose we have a solution $u$ of the nonautonomous problem in the
disease case such that $\partial_xu(t)$ is bounded in $X$. Then we
may write
\begin{eqnarray}
\label{full2} &&\partial_t u+\omega_\infty\partial_x u
+(\mu_0+\beta x) u-2\beta\int_x^\infty u(t,y)dy
=(\omega_\infty-\omega(t))\partial_xu,\\
&&u(0,x)=g(x),\quad u(t,0)=0,\quad t>0, x>0.\nonumber
\end{eqnarray}
Therefore we obtain the identity
$$u(t)=e^{-Lt}g +\int_0^t e^{-L(t-r)}(\omega_\infty-\omega(r))\partial_xu(r))dr,\quad t\geq0.$$
We know from Section 9 that $e^{-Lt}$ converges strongly in $X$ to
the ergodic projection $\cP$. On the other hand, the scalar
function $\omega(\cdot)-\omega_\infty$ belongs to $L_1(\R_+)$ by
assumption. This then implies
$$u(t)\rightarrow u_\infty\in R(\cP).$$
Thus we have convergence in $X$ to a unique element for all
nonnegative solutions with initial values in $D$. Since the
evolution operator associated with (\ref{full}) is bounded in $X$,
this convergence  extends to all initial values $u_0\in X$.

Returning now to the system \eqref{pde}, we may compute the limit
$u_\infty$. For this purpose recall that $U(t)=\int_{x_0}^\infty
u(t,x)dx\rightarrow U_\infty$ and $P(t)=\int_{x_0}^\infty
u(t,x)xdx\rightarrow P_\infty$. This implies
$$ u_\infty = \lim_{t\rightarrow\infty} \cP u(t)= \lim_{t\rightarrow\infty} [aU(t)+P(t)-x_0U(t)]e= [\mu U_\infty/\beta+P_\infty]e.$$
Note that $u_\infty$ is independent of the initial values $V_0$
and $u_0$.

This completes the proof of Theorem 1.2.

\bibliographystyle{amsnum}

\end{document}